# The messages we send after we listen

Yvonne Lai, Milton E. Mohr Associate Professor, University of Nebraska–Lincoln

Each time we teach, our students communicate to us. They talk, write, and draw. How do we take in their expressions of self and mathematics? How do we listen? In this column, I explore different ways we[1] as mathematics instructors may listen and how they contribute to a larger message of inclusion and induction into a mathematical community.

## When we listen to judge

One way we might listen is to judge. When we hear someone say, "All numbers are prime or composite," we might immediately think to ourselves: "That's wrong. They haven't said which numbers. And even if they meant positive integers, what about 1? Units are important in rings!" We judge for right and wrong.

Mathematics is sometimes said to be special because deductive logic is brittle to contradiction. If a statement has even one exception, then it is considered logically false. If a statement has any ambiguity, it is considered ill-defined. For a deductive argument to be solid, each statement must be entirely, unambiguously correct. We might listen for right and wrong because deductive beauty comes from knowing that your foundation is utterly, unbreakably solid.

Yet if we only ever listen to judge, and that is the only frame from which we respond, we might also send the message that participating in math is only about determining true and false. While understanding what it means to be true and false is important, this view of mathematics misses out on the sometimes exciting, sometimes frustrating, sometimes exulting messiness that is mathematical discovery.

There is a more insidious version of listening to judge as well. If we listen to judge whether a student is "worthy," we also send a message that some students are worthwhile and others are not. We need more students thinking that they can do mathematics, not fewer.

As mathematicians, we need to know how to assess the logical validity of a statement, and we want our students to know that too. At the same time, we need to be careful to use our judgment to help students rather than hurt students. We need to consider whether the way that we exercise our judgment supports students in building their mathematical experiences and selves.

## When we listen to understand

Every year when I teach Geometry for Secondary Teaching,[2] I ask prospective high school

---

[1] Throughout this essay, I use *we/our/us* as a stand-in for the community of those who teach mathematics, particularly at the tertiary level. The pronoun *I/my/me* is reserved for my personal actions and arguments.

[2] This course goes under many names depending on the institution, from Modern Geometry, to Geometry, to

mathematics teachers in the course to write definitions for *angle* and *angle measure*. Although I provide separate spaces on a worksheet to draft each, their definitions often blend together. I have the sense that the students see these concepts as so inextricably linked that they cannot untangle them to define them separately. Here are some sample drafts that groups have shared after comparing and building on individual group members' drafts:

- (Draft definition of *angle*) "The 2-d circular rotation from a vertex where we look at the difference from one ray to another."
- (Draft definition of *angle measure*) "The number that quantifies the circular path around the vertex of an angle."
- (Draft definition of *angle*) "The space between two rays that start at the same vertex, where the rays might overlap."

As you read these statements, you may be ready to judge. "That's not the definition of an angle! Why didn't they just say that an angle is the union of two rays with a common endpoint? Why doesn't Yvonne just give them the definition already?" You might be judging both my students' work and my pedagogical choices.

I want to suggest a different purpose for listening, even as we may still find ourselves judging. We can also listen to understand. We can set up situations where students can share their emerging ideas and practice communicating in writing what's in their heads. In other words, as mathematics educator Brent Davis (1997) wrote, we want to be "making sense of the sense [students] are making" (p. 365).[3] I ask the prospective high school math teachers to draft definitions of *angle* and *angle measure* because I want to know the images of these concepts that they bring with them.[4] Students are not blank slates. They come to us with all sorts of mathematical and emotional ideas of what mathematics and its concepts mean to them.

Some years ago, I went to a presentation on language and mathematics learning. The speaker said that the research in language acquisition suggests that the hardest technical vocabulary terms to learn are the ones that already have a different familiar meaning. In fact, the term *mean* is an example. This word can … mean … to define; an unkind action; and a number that is used to summarize values in a dataset such as an arithmetic, geometric, or harmonic mean. Children find the term *mean* harder to learn than, say, *standard deviation*, because there are so many other definitions and images that come up, other than an arithmetic mean, when hearing the word "mean." On the other hand, *standard deviation* likely holds no strong associations. Even though the term is a larger mouthful, its definition is easier to learn

---

Geometry for High School Teachers. Patricio Herbst and Amanda Brown, who lead the GeT: A Pencil group (https://getapencil.org/), have documented wide variation in these courses in topic and pedagogy. When I teach this course, I begin with a unit on foundational geometric arithmetic focused on vectors, angles, and betweenness. Then I turn to isometries of the Euclidean plane and their group theoretic properties, and then to congruence and similarity from a transformation perspective.

[3] Davis (1997) contrasted "interpretative" and "evaluative" listening. He characterized interpretative listening in terms of listening to make sense of students' own sense. He characterized evaluative listening in terms of a teacher who "already had a 'correct' answer in mind" (p. 360). Listening to judge is evaluative listening.

[4] In using the term "image," I allude to Vinner and Herchkowitz's (1980) notion of concept image, which includes all mental pictures of a concept in a student's mind, including associated processes and procedures, and how the pictures, processes, and procedures are connected.

because the space that the term occupies in the mind is as close to a blank slate as one can get.

I used to teach *angle* and *angle measure* by presenting a textbook definition. I found that my students clung so tightly to their own prior images that there was little room for any other. After some years of this, I wanted to surface their ideas. When I listened to their definitions to understand, I learned how the textbook definitions differed from their notions. What I learned from my students' drafts can be described in two recurring debates. First, students argue about whether an angle is the same thing as an angle measure. Some insist that the question, "What's the angle?" could legitimately be answered "60 degrees" or another number of degrees. These students think of the angle as the measure, rather than as a physical object in and of itself. Second, students argue about whether the interior of an angle is part of the angle or not. Some students think an angle is just the rays. Other students think an angle is the rays and a region inside. Still other students think that the angle is only the interior. These debates came up because students compared and contrasted their conceptions with each other and the textbook. With a clear sense of the prior images held by students, I was able to say directly how their definitions and the textbook definition differed. We then talked about the utility of each, and a rationale for the textbook's choice.

I want to understand prospective high school teachers' personal concepts of *angle* and *angle measure* because I want future teachers to connect angles in geometry to angles in precalculus and trigonometry. To do so, I want to have a conversation about why the "inside" of an angle might be something included in the definition of an angle. Related to this, I also want to know whether they see angle measure as only non-negative quantities or whether they see angle measure as signed. By taking angle measure as signed rather than non-negative, geometry connects to trigonometry, and one can also define rotations precisely. Moreover, properties of compositions of reflections and rotations are more elegantly expressed with signed angle measures. The composition of two rotations has an angle measure that is simply the addition of their individual measures, rather than an addition with separate cases for clockwise and counterclockwise.[5]

In general, I believe that we should strive to know when we go against student intuitions, so that we can work with our students to develop new intuitions. The only way we can know whether we go against student intuitions is by listening to understand. The more we can help students develop new intuitions, the more learning possibilities we open.

## When we listen to create

Here are some draft definitions of "triangle" written by students in the geometry class:[6]

- "A [closed] 2-dimensional shape with 3 vertices whose interior angle measures add up to 180 degrees." (This group inserted "closed" after reading other groups' definitions.)

---

[5] This fact is beautifully true even if the rotations do not share the same center. For instance, if you have two points $X$ and $Y$ in the plane, and you compose a rotation of –30 degrees about $X$ with a rotation of +60 degrees about $Y$, there is some point $Z$ such that the result of the composition is a rotation of +30 degrees about $Z$. (The point $Z$ is different depending on the order of composition.) If $X = Y$, then $Z = X = Y$.

[6] I wrote about this activity in the September-October 2022 Education Column (Lai, 2022). Here I share some more student ideas that came out of this activity.

- "A closed shape with three straight sides."
- "A closed 2-dimensional shape made of exactly three straight edges that are joined endpoint to endpoint."

As I looked around the room that day at students' definitions, many were drawn to using the term "closed." As I listened to understand, I asked them questions such as, "When you say 'closed,' what do you mean?", "What's a drawing of something else that's 'closed'?", "What's a drawing of something that's not 'closed'?", "What made you think of using 'closed' in your definition?" Overall, students were not sure why, but they felt convinced that a correct definition needed that word.

I saw an opportunity to create a definition together. Although I had not planned to do so, I switched up the groups and asked the new groups to define "closed" and to draw examples of what "closed" meant. Each group claimed a whiteboard. Our discussion focused on planar curves. The students pointed to existing figures in the room and drew new ones, as in Figure 1. They began to bar some from being "closed" as though they were monstrous to the concept's meaning.[7] Using their suggestions and critiques, I drafted and revised a proposed definition. The class came to consensus on a definition based on being able to trace the entire shape and return to where you started, without any backtracking, and allowing for a finite number of isolated point overlaps in the tracing. Some students, who were taking a graph theory class concurrently, connected our conversation to circuits in graphs.

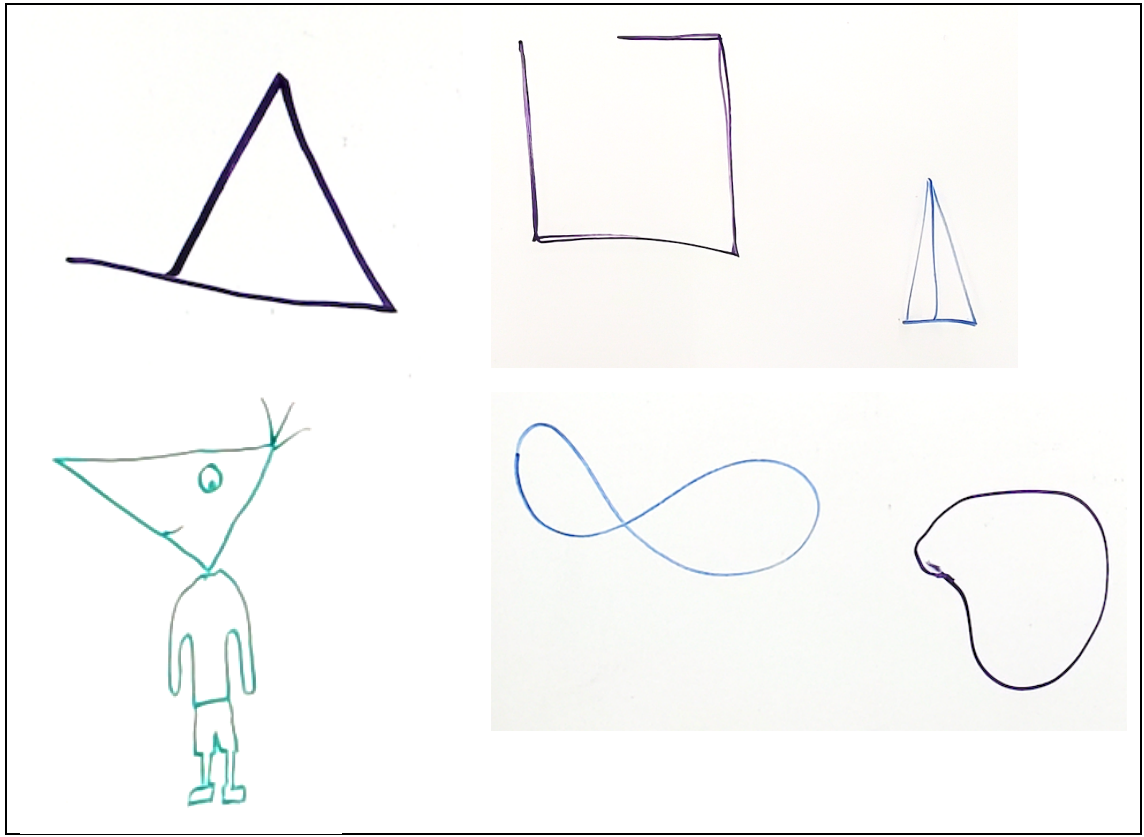

---

[7] Imre Lakatos (1963) described the process of excluding examples as "monster barring."

**Figure 1.** *Drawings that students made along the way to creating a definition of closed.*

Listening to create means listening to spot opportunities to build mathematics together.[8] That mathematics may take the form of a definition, conjecture, proof, mathematical model, or another way to communicate and advance the community's knowledge.

## When we listen through desire, and learning to listen past desire

Before closing this essay, I will issue a caution. Sometimes when we think we are listening to understand or create, we may be listening through desire. Once, when a student asked me about "0" in a ring theory class, I launched into a spontaneous lecture on aspects of the zero-product property. The student gaped. It turned out that due to using an old whiteboard marker, what I intended as an alpha, they read as a zero. They were confused why a zero would show up in the calculation. It is an enduring problem of teaching to learn to listen past desire.

## Who gets to participate, and what does participation look like?

How we listen answers the question: *Who gets to participate in mathematics, and what does participation look like?* When we listen to judge for true, false, or well-definedness, we send the message that precision is important in doing mathematics. At the same time, if we only listen to judge in this way, we send the message that participating in mathematics can only be done with precise language. Those who falter with imprecise language are not participating. Yet inquiry in mathematics, whether we seek a proof or a mathematical model, relies on a place for inchoate ideas that are refined into more precise communication. Without space for imprecision, the precision would mean little.

When we listen to understand, we send the message that students' concept images matter to where we go next. Participation in mathematics involves surfacing, comparing, and contrasting conceptions of mathematical ideas. In doing so, we can analyze choices made in mathematical definitions and their implications for later results. When students discuss the choices made in mathematical definitions and assumptions, they are participating in defining the terms of justification in the class community. By listening to understand, we honor students' ideas and set them up to create new ideas.

When we listen to create, we draw on the skills we gain from listening for precision and listening to understand. We also go beyond these skills to show that the agenda of mathematics is determined not by the instructor's fiat, but by a community. We send the message that participation in mathematics means shaping fuzzy ideas into sharper questions and claims, understanding where each other come from, and seeing what we can build together.

## Acknowledgements

[8] "Listening to create" is my interpretation of Yackel et al.'s (2003) "generative listening" and Davis's (1997) "hermeneutic listening." Yackel et al. (2003) defined generative listening as improvising a revised lesson trajectory based on student contributions so as to "generate or transform one's own mathematical understanding and … generate a new space of instructional activities" (p. 117). Davis (1997) defined hermeneutic listening as negotiated and participatory, involving a "willingness to interrogate the taken for granted" (pp. 369-370).

*The framework of listening to judge, understand, and create comes from a presentation by Gail Burrill at the 2019 Joint Mathematics Meeting, where she gave references to the frameworks of Brent Davis (1997) and of Erma Yackel and colleagues (2003). She spoke at a panel on Listening and Responding to Students in K-16, held by the Special Interest Group of the MAA on Mathematical Knowledge for Teaching (SIG-MKT). Ever since hearing her presentation, I have thought about how these modes of listening impact the messages we send to students about what it means to participate in mathematics. This essay is an expansion of these thoughts, which I first had an opportunity to discuss with others at a presentation in the National Museum of Mathematics' virtual QED series. I am grateful to Gail Burrill for introducing me to the framework, to Tim Chartier for the invitation to speak, and to Christopher Danielson and Francis Su for encouraging me to explore these ideas in writing.*

**Keywords:** mathematics teaching practice, proof, listening, inclusion